%% file: Loevbak.tex
\pdfoutput=1
\documentclass[graybox]{svmult}

\usepackage{type1cm}        
\usepackage{makeidx}         
\usepackage{graphicx}        
\usepackage{multicol}        
\usepackage[bottom]{footmisc}

\usepackage{newtxtext}       %
\usepackage[varvw]{newtxmath}       

\makeindex             

\usepackage{amsmath}
\usepackage{autonum}
\usepackage{sourcecodepro}
\usepackage{listings}

\usepackage{dirtytalk}
\usepackage{tikz}
\usepackage{pgfplots}
\pgfplotsset{compat=1.17}
\usepgfplotslibrary{fillbetween}
\usepackage{pgfplotstable}
\usepackage{xcolor}
\usepackage[caption=false]{subfig}

\begin{document}

\title*{Reversible random number generation for adjoint Monte Carlo simulation of the heat equation}
\titlerunning{Reversible RNG for adjoint MC simulation of the heat equation}
\author{Emil Løvbak \and Frédéric Blondeel \and Adam Lee \and Lander Vanroye \and Andreas Van Barel \and Giovanni Samaey}
\authorrunning{Emil Løvbak et al.}
\institute{Emil Løvbak \and Frédéric Blondeel \and Adam Lee \and Andreas Van Barel  \and Giovanni Samaey \at KU Leuven, Department of Computer Science, Celestijnenlaan 200A box 2402, 3001 Leuven \\ \email{emil.loevbak@kuleuven.be $\cdot$ fredericblondeel@hotmail.com  $\cdot$ adamrlee@pm.me \\ andreas.vanbarel@gmail.com $\cdot$ giovanni.samaey@kuleuven.be}
\and Lander Vanroye \at KU Leuven, Department of Mechanical Engineering, Celestijnenlaan 300 box 2420, 3001 Leuven \\ \email{lander.vanroye@kuleuven.be}}
%
%
\maketitle

\abstract*{In PDE-constrained optimization, one aims to find design parameters that minimize some objective, subject to the satisfaction of a partial differential equation. A major challenge is computing gradients of the objective to the design parameters, as applying the chain rule requires computing the Jacobian of the design parameters to the PDE's state. The adjoint method avoids this Jacobian by computing partial derivatives of a Lagrangian. Evaluating these derivatives requires the solution of a second PDE with the adjoint differential operator to the constraint, resulting in a backwards-in-time simulation.
Particle-based Monte Carlo solvers are often used to compute the solution to high-dimensional PDEs. However, such solvers have the drawback of introducing noise to the computed results, thus requiring stochastic optimization methods. To guarantee convergence in this setting, both the constraint and adjoint Monte Carlo simulations should simulate the same particle trajectories. For large simulations, storing full paths from the constraint equation for re-use in the adjoint equation becomes infeasible due to memory limitations. In this paper, we provide a reversible extension to the family of permuted congruential pseudorandom number generators (PCG). We then use such a generator to recompute these time-reversed paths for the heat equation, avoiding these memory issues.
}

\abstract{
In PDE-constrained optimization, one aims to find design parameters that minimize some objective, subject to the satisfaction of a partial differential equation. A major challenge is computing gradients of the objective to the design parameters, as applying the chain rule requires computing the Jacobian of the design parameters to the PDE's state. The adjoint method avoids this Jacobian by computing partial derivatives of a Lagrangian. Evaluating these derivatives requires the solution of a second PDE with the adjoint differential operator to the constraint, resulting in a backwards-in-time simulation.\newline\indent
Particle-based Monte Carlo solvers are often used to compute the solution to high-dimensional PDEs. However, such solvers have the drawback of introducing noise to the computed results, thus requiring stochastic optimization methods. To guarantee convergence in this setting, both the constraint and adjoint Monte Carlo simulations should simulate the same particle trajectories. For large simulations, storing full paths from the constraint equation for re-use in the adjoint equation becomes infeasible due to memory limitations. In this paper, we provide a reversible extension to the family of permuted congruential pseudorandom number generators (PCG). We then use such a generator to recompute these time-reversed paths for the heat equation, avoiding these memory issues.
}

\section{Introduction}

Partial differential equations (PDEs) are an indispensable tool for modeling physics in many engineering fields. When solving a design problem, e.g., building the cheapest, most efficient component for a given purpose, the design must satisfy the problem's physics, requiring a numerical solver for the PDE. To improve the design's objective function, one can compute its gradient to the design parameters using the discrete adjoint approach. If the PDE solver is based on a particle simulation, then the adjoint solver traces the same particle trajectories in a time reversed fashion. Storing all of these paths can incur unacceptably high memory costs. We therefore propose a novel approach to recompute these paths using a reversible random number generator.

Mathematically, we consider a PDE-constrained optimization problem of the form 
\begin{equation}
	\label{LoBlLe:eq:constrainedoptimization}
	\min_{u(x)} \mathcal{J}(y(x,t), u(x)), \quad \text{subject to} \quad  \mathcal{B}(y(x,t); u(x)) = 0.
\end{equation}
Here, $\mathcal{B}(\cdot; u(x))$ is a differential operator, parameterized by a time-independent control $u(x)$, applied to a candidate solution $y(x,t)$, a function of space $x \in \mathcal{D}_x \subseteq \mathbb{R}^d$ and time $t \in \mathbb{R}_{\geq 0}$. The objective $\mathcal{J}(y(x,t), u(x))$, is a scalar function to be minimized, typically involving integration over space and time. For all but the most trivial of applications, solving the PDE requires a discretization. When discretizing the PDE  $\mathcal{B}(y(x,t); u(x))$ and objective$\mathcal{J}(y(x,t), u(x))$, one gets a discrete equivalent to \eqref{LoBlLe:eq:constrainedoptimization},
\begin{equation}
	\label{LoBlLe:eq:discreteconstrainedoptimization}
	\min_{\hat{u}} \hat{\mathcal{J}}(\hat{y}, \hat{u}), \quad \text{subject to} \quad  \hat{\mathcal{B}}(\hat{y}; \hat{u}) = 0,
\end{equation}
where $\hat{y} \in \mathbb{R}^\alpha$, $\hat{u} \in \mathbb{R}^\beta$, $\hat{\mathcal{J}} \colon \mathbb{R}^\alpha \times \mathbb{R}^\beta \to \mathbb{R}$ and $\hat{\mathcal{B}}(\cdot; \hat{u}) \colon \mathbb{R}^\alpha \to \mathbb{R}^\gamma$ are discretizations of the quantities in \eqref{LoBlLe:eq:constrainedoptimization} and the integer values $\alpha$, $\beta$ and $\gamma$ are discretization dependent.

In this work, we solve the PDE with particle-based Monte Carlo, i.e., we simulate sample trajectories of particles whose distribution density corresponds with a rescaling of the PDE's solution. The main advantage of Monte Carlo methods lies in their ability to simulate high-dimensional PDEs, without storing their state on a high-dimensional grid. These methods, however, come with the drawback of introducing noise on the computed results, the variance of which scales as $P^{-1}$, with $P$ the number of simulated particle trajectories. Stochastic optimization routines are therefore needed to solve \eqref{LoBlLe:eq:discreteconstrainedoptimization}. While this work is motivated by its future applicability in the simulation of kinetic models in neutral particle codes for fusion reactor design as used in~\cite{Feng1997, Reiter2005}, we consider a simplified setting here, based on the heat equation.

Solving \eqref{LoBlLe:eq:discreteconstrainedoptimization} with a gradient-based optimization method, e.g., stochastic gradient descent, requires evaluating both $\hat{\mathcal{J}}(\hat{y}^\prime,\hat{u}^\prime)$ and $\frac{\text{d} \hat{\mathcal{J}}}{\text{d} \hat{u}} (\hat{y}^\prime, \hat{u}^\prime)$ for a proposed control $\hat{u}^\prime$ and corresponding state $\hat{y}^\prime$, with $ \hat{\mathcal{B}}(\hat{y}^\prime; \hat{u}^\prime) = 0$. As the numerical solver $\hat{\mathcal{B}}(\cdot; \hat{u}^\prime)$ defines an implicit function $\hat{y}^\prime = \hat{y}(\hat{u}^\prime)$, one may try to compute the gradient as
\begin{equation}
	\label{LoBlLe:eq:chainrule}
	\frac{\text{d} \hat{\mathcal{J}}}{\text{d} \hat{u}} (\hat{y}(\hat{u}^\prime), \hat{u}^\prime) = \frac{\partial \hat{\mathcal{J}}}{\partial \hat{u}}(\hat{y}^\prime, \hat{u}^\prime) + \frac{\partial \hat{\mathcal{J}}}{\partial \hat{y}}(\hat{y}^\prime, \hat{u}^\prime) \frac{\text{d} \hat{y}}{\text{d} \hat{u}}(\hat{u}^\prime).
\end{equation}
Evaluating the Jacobian $\frac{\text{d} \hat{y}}{\text{d} \hat{u}}(\hat{u}) \in \mathbb{R}^{\alpha \times \beta}$ is often infeasible, given its dimentionality and the implicit definition of $\hat{y}(\hat{u})$. We thus consider three alternatives to compute $\frac{\text{d} \hat{\mathcal{J}}}{\text{d} \hat{u}}$: finite differences, algorithmic differentiation and the adjoint-based approach.

The finite-difference approach is the only one of the three which is a true black-box approach, i.e., not needing knowledge of the solver internals. It is also simple to implement. For each element of $\frac{\text{d} \hat{\mathcal{J}}}{\text{d} \hat{u}}(\hat{y}^\prime, \hat{u}^\prime)$, one simply subtracts two solver outputs from each other and performs a division. Its main drawbacks are that it scales poorly with the dimension of $\hat{u}$, as the PDE needs to be solved $\beta$ times in addition to the solve for computing $\hat{y}^\prime$, and that it typically only achieves a relative precision in the order of $\sqrt{\epsilon_\text{mach}}$ due to cancellation errors, with $\epsilon_\text{mach}$ the machine precision.

Algorithmic differentiation~\cite{Gebremedhin2020}, while not fully black-box, takes the existing solver code and automatically generates code for calculating its derivative. Conceptually, this approach relies on repeated application of the chain rule. In our case, so-called backward mode is the most advantageous as $\beta \gg 1$. This approach computes $\frac{\text{d} \hat{\mathcal{J}}}{\text{d} \hat{u}}(\hat{y}^\prime, \hat{u}^\prime)$ to a relative precision in the order of $\epsilon_\text{mach}$, with a theoretical cost upper-bound of four times that of evaluating $\hat{\mathcal{J}}(\hat{y}^\prime, \hat{u}^\prime)$~\cite{Gebremedhin2020}. However, in the context of Monte Carlo simulations, this approach requires accessing particle states in a time-reversed fashion. Storing these paths in memory is infeasible for all but toy problems. Therefore, one soon needs to apply checkpointing, i.e., storing the simulation state at intermediate moments in time from which partial paths can be recomputed as they are needed~\cite{Naumann2018}. This recomputation has at least the cost of the initial constraint simulation and the checkpointing itself often incurs a non-negligible overhead.

The adjoint-based approach~\cite{Giles2000} uses a solver for $\hat{\mathcal{B}}^\ast(\cdot; \hat{y}^\prime, \hat{u}^\prime)$, the adjoint operator to $\hat{\mathcal{B}}(\cdot; \hat{u}^\prime)$, derived at $\hat{y}^\prime$. In this approach, we compute the gradient $\frac{\text{d} \hat{\mathcal{J}}}{\text{d} \hat{u}} (\hat{y}^\prime, \hat{u}^\prime)$, through a vector of Lagrange multipliers $\hat{\lambda}$. The main advantage of this approach is the reduced computational cost. The adjoint PDE for computing $\hat{\lambda}$ is typically similar to the constraint PDE, so the gradient can be computed at the cost of one additional PDE solve and a small number of matrix-vector products. The main drawback of the adjoint-based approach is the additional work needed in deriving and implementing the adjoint solver. However, as the constraint and adjoint PDE are often similar, existing code can be reused with minimal modifications to develop such a solver.

In the adjoint-based approach one can apply two strategies~\cite{Caflisch2021}. In the first strategy, called the continuous adjoints or optimize then discretize (OTD), one first derives the adjoint PDE and discretizes each PDE independently. When using Monte Carlo simulation, this strategy has two issues. Often, the adjoint PDE depends on the solution of the constraint PDE, requiring the construction of a solution on the high-dimensional grid that we wish to avoid. Discretizing both PDEs seperately also decouples the resulting stochastic simulations, hampering the convergence of stochastic gradient descent, see e.g.~\cite[Thm. 4.8]{Bottou2018}. We therefore use discrete adoints or discretize then optimize (DTO), where we directly derive an adjoint solver from the constraint solver. This derivation results in the same particle paths being used in both discretizations, except for the adjoint simulation running backwards in time.

Note that DTO avoids computing a high-dimensional solution, while introducing a need to store the Monte Carlo paths used in solving the constraint PDE, for use in the adjoint simulation. Note that this is the same issue as was present in backward-mode algorithmic differentiation, which has conceptual similarities to the adjoint-based approach. While checkpointing can be used to avoid the memory issue, it has the additional drawback that it introduces extra complexity to the adjoint simulation. 

The goal of this work is twofold: to present a DTO adjoint approach for Monte Carlo simulations, and to do away with checkpointing by re-computing paths time-reversed during the adjoint simulation. To compute time-reversed paths, we extend the family of permuted congruential pseudorandom number generators (PCG)~\cite{ONeill2014} to make them reversible, i.e., one can request the next and previous random value in the sequence at the same computational cost. While a similarly reversible generator was developed in~\cite{Yoginath2018}, based on Multiple Recursive Generators~\cite{LEcuyer2000}, we opt for our own implementation based on PCG due to its performance and statistical properties.

The remainder of this paper is structured as follows. In Section~\ref{LoBlLe:sec:adjointMC}, we introduce the discrete adjoint approach in the context of Monte Carlo simulations. We then introduce a 1D test problem, based on the heat equation, on which we demonstrate our reversible approach. In Section~\ref{LoBlLe:sec:rrng}, we introduce the PCG family of random number generators and show how we adapt the generators to make them reversible. We then demonstrate the speedup and memory gains provided by the reversible approach in Section~\ref{LoBlLe:sec:experiments}. Finally in Section~\ref{LoBlLe:sec:conclusions}, we draw our conclusions and outline our future work.

\section{Adjoint Monte Carlo for particle simulations}
\label{LoBlLe:sec:adjointMC}

In the adjoint-based optimization approach, we replace the constrained discrete problem \eqref{LoBlLe:eq:discreteconstrainedoptimization} by one of finding stationary points of the discrete Lagrangian
\begin{equation}
	\label{LoBlLe:eq:lagrangian}
	\hat{\mathcal{L}}(\hat{y}, \hat{u}, \hat{\lambda}) = \hat{\mathcal{J}}(\hat{y}, \hat{u})+ \hat{\lambda}^\top \hat{\mathcal{B}}(\hat{y}; \hat{u}),
\end{equation}
with $\hat{\lambda} \in \mathbb{R}^\gamma$ a vector of Lagrange multipliers. We now find these fixed points by setting the partial derivatives in $\hat{\lambda}$, $\hat{y}$ and $\hat{u}$ to zero:

\begin{align}
	\frac{\partial\hat{\mathcal{L}}}{\partial\hat{\lambda}}^\top\!\!\!(\hat{y}, \hat{u},\hat{\lambda}) &= \hat{\mathcal{B}}(\hat{y}; \hat{u}) = 0, && \quad \text{Constraint equation}; \label{LoBlLe:eq:state}\\
	\frac{\partial\hat{\mathcal{L}}}{\partial\hat{y}}^\top\!\!\!(\hat{y}, \hat{u},\hat{\lambda}) &= \frac{\partial\hat{\mathcal{J}}}{\partial\hat{y}}^\top\!\!\!(\hat{y}, \hat{u}) + \frac{\partial\hat{\mathcal{B}}}{\partial\hat{y}}^\top\!\!\!(\hat{y};\hat{u})\hat{\lambda} = 0, && \quad \text{Adjoint equation}; \label{LoBlLe:eq:adjoint}\\
	\frac{\partial\hat{\mathcal{L}}}{\partial\hat{u}}^\top\!\!\!(\hat{y}, \hat{u},\hat{\lambda}) &= \frac{\partial\hat{\mathcal{J}}}{\partial\hat{u}}^\top\!\!\!(\hat{y}, \hat{u}) + \frac{\partial\hat{\mathcal{B}}}{\partial\hat{u}}^\top\!\!\!(\hat{y};\hat{u})\hat{\lambda}   =0, && \quad \text{Design equation}. \label{LoBlLe:eq:design}
\end{align}
From \eqref{LoBlLe:eq:state}, we see that the constraining PDE is satisfied for any stationary point of the Lagrangian. Therefore, it holds that, for any $\hat{y}^\prime$, $\hat{u}^\prime$ and $\hat{\lambda}^\prime$ solving the system \eqref{LoBlLe:eq:state}--\eqref{LoBlLe:eq:design}, $\hat{\mathcal{L}}(\hat{y}^\prime, \hat{u}^\prime, \hat{\lambda}^\prime) = \hat{\mathcal{J}}(\hat{y}^\prime, \hat{u}^\prime)$. A fixed point of the Lagrangian is thus also a fixed point of the objective function and a possible solution to the constrained problem \eqref{LoBlLe:eq:discreteconstrainedoptimization}.

For Monte Carlo, the simulation state $\hat{y}$ consists of an ensemble of $P$ particles
\begin{equation}
	\label{LoBlLe:eq:ensemble}
	\left\{ Y_{p,\tau} \right\}_{p=1}^P,
\end{equation}
with $Y_{p,\tau}$ a vector valued quantity defined at each moment in time $\tau = 0,\dots, T$ for each particle with index $p=1, \dots, P$. Assuming independent particles, their time discretized dynamics are described by algebraic expressions 
\begin{equation}
\label{LoBlLe:eq:algebraicform}
B_{p,\tau+1}(Y_{p,\tau+1}, Y_{p,\tau}; \hat{u}) = 0, \quad \tau = 0,\dots, T-1,
\end{equation}
relating particle states at subsequent time steps.We structure the simulation state as $\hat{y} = \begin{bmatrix} \begin{bmatrix} Y_{0,1}^\top & \cdots & Y_{0,T}^\top \end{bmatrix} & \cdots & \begin{bmatrix} Y_{P,1}^\top & \cdots & Y_{P,T}^\top \end{bmatrix} \end{bmatrix}^\top$, with the system $ \hat{\mathcal{B}}(\hat{y}; \hat{u})=0$ being given by the aggregation of the expressions \eqref{LoBlLe:eq:algebraicform}. We now derive a linearization
\begin{equation}
	\label{LoBlLe:eq:linearization}
	\hat{\mathcal{B}}(\hat{y}; \hat{u}^\prime) \approx \frac{\partial\hat{\mathcal{B}}}{\partial\hat{y}}(\hat{y}^\prime, \hat{u}^\prime) (\hat{y}-\hat{y}^\prime)
\end{equation}
 of the system assuming we have computed a solution $\hat{y}^\prime$ for a given $\hat{u}^\prime$ so that $\hat{\mathcal{B}}(\hat{y}^\prime; \hat{u}^\prime)=0$. We now assume piecewise linearity of $\hat{\mathcal{B}}(\hat{y}; \hat{u})$ in $\hat{y}$ making \eqref{LoBlLe:eq:linearization} exact for a finite region of values $\hat{y}$ around $\hat{y}^\prime$. This approach is also applicable to the non-linear case, under suitable combination with an iterative solver for \eqref{LoBlLe:eq:state} and \eqref{LoBlLe:eq:adjoint}.

In \eqref{LoBlLe:eq:linearization}, we see that a Monte Carlo simulation, within a region of the solution space, can be written as a matrix-vector product between a Jacobian and a vector of offset particle states. We also observe that the Jacobian $\frac{\partial\hat{\mathcal{B}}}{\partial\hat{y}}(\hat{y},\hat{u})$ is blockwise lower-triangular with a single off-diagonal band of blocks, which is solved for $\hat{y}$ using forward substitution. We now observe that the same Jacobian appears, transposed, in \eqref{LoBlLe:eq:adjoint}. We can thus solve \eqref{LoBlLe:eq:adjoint} for $\hat{\lambda}$ with a similar Monte Carlo simulation to \eqref{LoBlLe:eq:algebraicform}. The transposed Jacobian is upper-triangular, meaning the adjoint Monte Carlo simulation corresponds with backward substitution, and thus runs backward in time.

In Section~\ref{LoBlLe:sec:1Dheat}, we present a test problem, using the heat equation. In Section~\ref{LoBlLe:sec:adjointchallenges}, we discuss adjoint simulation time reversal, motivating our reversible approach.

\subsection{Cooling a 1D rod}
\label{LoBlLe:sec:1Dheat}

We consider the problem of cooling a 1D rod with temperature $\theta(x,t)$ with $x \in [0,L]$ and $t \in \mathbb{R}_{\geq 0}$, cooled by a fluid with a rate $u(x)$, which we write formally as
\begin{align}
	\min_{u(x)} \; \mathcal{J}(\theta(x,t), u(x,t)) &= \int_0^\infty \int_0^L \frac{1}{2}\theta(x,t)^2 \text{d} x \, \text{d} t + \nu \int_0^L \frac{1}{2}u(x)^2 \text{d} x, \label{LoBlLe:eq:heatobjective}\\
	 \text{subject to} \;\mathcal{B}(\theta(x,t); u(x,t)) &= \frac{\partial}{\partial t} \theta(x,t) - \frac{\partial^2}{\partial x^2} \theta(x,t) + u(x)\theta(x,t) = 0,\label{LoBlLe:eq:heat}\\
	  \theta(x,0) = \theta_0(x), &\quad \theta(0,t) = \theta(L,t), \label{LoBlLe:eq:heatinitial}
\end{align}
where $\nu \in \mathbb{R}_{\geq 0}$ is a regularization parameter. For simplicity, we consider $\theta(x,t)$ and $u(x)$ to be dimensionless.

Given a spatial discretization with cell width $\Delta x$, we perform the approximation
\begin{equation}
\label{LoBlLe:eq:heatdiscreteobjective}
\mathcal{J}(\theta(x,t), u(x,t)) \approx \hat{\mathcal{J}}(\hat{\theta}, \hat{u}) = \Delta t \sideset{}{''} \sum_{\tau=0}^{T} \Delta x \frac{1}{2}\hat{\theta}_\tau^\top \hat{\theta}_\tau + \nu  \Delta x \frac{1}{2}\hat{u}^\top \hat{u},
\end{equation}
where we introduce $\hat{\theta} \in \mathbb{R}^{(T+1) \times N}$ and $\hat{u} \in \mathbb{R}^N$ as vector valued quantities indexed over time $\tau = 0,\dots, T$ and space $n=1,\dots,N$ so that $N \Delta x = L$, and $T\Delta t$ is sufficiently large that \eqref{LoBlLe:eq:heat} approximates a steady-state solution. Note that $\hat{\theta}_{\tau,n}$ denotes the solution value in the cell $(n-1)\Delta x \leq x < n\Delta x$ at time $t=\tau\Delta t$. The notation $\sum\nolimits''$ indicates halving of the first and last terms due to integration with the trapezoid rule.

At the Monte Carlo level, we consider an ensemble of $P$ \say{heat-transfer} particles
\begin{equation}
	\label{LoBlLe:eq:heatensemble}
	\left\{ Y_{p,\tau} \right\}_{p=1}^P = \left\{ [X_{p,\tau}, W_{p,\tau}]^\top \right\}_{p=1}^P,
\end{equation}
determined by their position $X_{p,\tau}$ and a weight $W_{p,\tau}$ at each time step $\tau$. Given an ensemble \eqref{LoBlLe:eq:heatensemble}, one can apply a binning approach to compute a stochastic quantity $\hat{\Theta}$ for which $\mathbb{E}\left[\hat{\Theta}\right] = \hat{\theta}$, where for each $\tau$, $n$
\begin{equation}
	\hat{\Theta}_{\tau,n} = \sum_{p=1}^P \frac{1}{\Delta x} \mathcal{I}_n\left( X_{p,\tau} \right) W_{p,\tau},
\quad
	\text{with} \quad \mathcal{I}_n(x) = \begin{cases}
		1, & \text{if} \, (n-1)\Delta x \leq x < n\Delta x;\\
		0, & \text{otherwise}.
	\end{cases}
\end{equation}
By substituting the stochastic estimate $\hat{\Theta}$ for $\hat{\theta}$ in \eqref{LoBlLe:eq:heatdiscreteobjective}, one then gets a stochastic estimate for $\hat{\mathcal{J}}$ for a given ensemble of particles \eqref{LoBlLe:eq:ensemble}.

To initialize the simulation, we sample each particle's initial position $X_{p,0}$ from the density function $\frac{\theta_0(x)}{\int_0^L \theta_0(x) \text{d}x}$. We then select initial weights $W_{p,0}$, such that $\hat{\Theta}_0=\hat{\theta}_0$. We then iterate over a particle scheme consisting of two parts for each time step:
\begin{itemize}
	\item \textbf{Diffusion step.} Update the particle's position following a Brownian motion by sampling a normally distributed value $\xi_{p,\tau}$, i.e.,
	\begin{equation}
		\label{LoBlLe:eq:diffusionstep}
		X_{p,\tau+1} = X_{p,\tau} + \sqrt{2\Delta t}\xi_{p,\tau} \quad \xi_{p,\tau} \sim \mathcal{N}(0,1).
	\end{equation}
	\item \textbf{Reweighting step.} Reduce the particle's weight with an exponential decay, i.e.,
	\begin{equation}
		\label{LoBlLe:eq:weightupdate}
		W_{p,\tau+1} = W_{p,\tau} \exp\left(-\Delta t\,\hat{u}\left(X_{p,\tau+1}\right)\!\right),
	\end{equation}
	where $\hat{u}(x)$ denotes the value $\hat{u}_n$ with $n$ chosen so that $\mathcal{I}_n(x)$ is non-zero.
\end{itemize}
Each particle trajectory $p$ is thus fully determined by its initial position $X_{p,0}$ and its Brownian increments $\xi_{p,\tau}$.

Rewriting the time-stepping in the form \eqref{LoBlLe:eq:algebraicform}, for a given $\hat{u}^\prime$ gives
\begin{equation}
	\label{LoBlLe:eq:heatrecursion}
	B_{p,\tau+1}(Y_{p,\tau+1}, Y_{p,\tau}; \hat{u}^\prime) = \begin{bmatrix}
		X_{p,\tau+1} - X_{p,\tau} - \sqrt{2\Delta t}\xi_{p,\tau} \\
		W_{p,\tau+1} -W_{p,\tau} \exp\left(-\Delta t\,\hat{u}^\prime\!\!\left(X_{p,\tau+1}\right)\!\right)
	\end{bmatrix}, \; \tau = 0, \dots, T-1.
\end{equation}
The Jacobian $\frac{\partial\hat{\mathcal{B}}}{\partial\hat{y}}(\hat{y}^\prime,\hat{u}^\prime)$ is then block-diagonal with blocks $p = 1,\dots, P$ given by
\begin{equation}
	\frac{\partial\hat{\mathcal{B}}_p}{\partial\hat{y}_p}(\hat{y}_p^\prime,\hat{u}^\prime) = \begin{bmatrix}
		\frac{\partial B_{p,1}}{\partial Y_{p,1}} \\
		 \frac{\partial B_{p,2}}{\partial Y_{p,1}} & \frac{\partial B_{p,2}}{\partial Y_{p,2}} \\
		 & \ddots & \ddots \\
		 &	& \frac{\partial B_{p,T}}{\partial Y_{p,T-1}} & \frac{\partial B_{p,T}}{\partial Y_{p,T}} \\
	\end{bmatrix}\!,
\,	\text{with} \,
	\frac{\partial B_{p,\tau+1}}{\partial Y_{p,\tau+1}}\left(Y_{p,\tau+1}, Y_{p,\tau}; \hat{u}^\prime\right) = \begin{bmatrix}
		1 & 0 \\
		0 & 1
	\end{bmatrix}
\end{equation}
\begin{equation}
	\text{and} \quad
	\frac{\partial B_{p,\tau+1}}{\partial Y_{p,\tau}}\left( Y_{p,\tau+1}, Y_{p,\tau}; \hat{u}^\prime\right) = \begin{bmatrix}
		-1 & 0 \\
		0 &  -\exp\left(-\Delta t\,\hat{u}^\prime\!\!\left(X_{p,\tau+1}\right)\!\right) 
	\end{bmatrix},
\end{equation}
where the bottom-left element of $\frac{\partial B_{p,\tau+1}}{\partial Y_{p,\tau+1}}$ is zero by the piecewise constantness of $\hat{u}$.

From \eqref{LoBlLe:eq:adjoint}, we derive a final condition at $\tau=T$ for each particle $p$
\begin{equation}
	\label{LoBlLe:eq:heatadjointend}
	\hat{\lambda}_{p,T} = \begin{bmatrix}
		X^\ast_{p,T} \\ W^\ast_{p,T}
	\end{bmatrix} = - \frac{\partial \hat{\mathcal{J}}}{\partial Y_{p,T}}^{\!\!\!\!\top} = - \begin{bmatrix} 0 \\ \Delta t \Theta \left( X_{p,T}\right) \end{bmatrix},
\end{equation} 
introducing $X^\ast_{p,T}$ and $W^\ast_{p,T}$ as adjoint variables, and a recursion for $\tau = T-1,\dots,1$
\begin{equation}
	\label{LoBlLe:eq:heatadjoint}
	\hat{\lambda}_{p,\tau} = - \frac{\partial B_{p,\tau+1}}{\partial Y_{p,\tau}}^{\!\!\!\top} \hat{\lambda}_{p,\tau+1}  - \frac{\partial \hat{\mathcal{J}}}{\partial Y_{p,\tau}}^{\!\!\!\!\top} = \begin{bmatrix} X^\ast_{p,\tau+1} \\ \exp\left(-\Delta t\,\hat{u}^\prime\!\!\left(X_{p,\tau+1}\right)\!\right) W^*_{p,\tau+1} \end{bmatrix} - \begin{bmatrix} 0 \\ \Delta t \Theta(X_{p,\tau}) \end{bmatrix} .
\end{equation}
Observe that $\forall p,\tau : X^\ast_{p,\tau} \equiv 0$ meaning \eqref{LoBlLe:eq:heatadjoint} defines a Monte Carlo simulation
\begin{equation}
	\label{LoBlLe:eq:heatreversemc}
	W^\ast_{p,\tau} =  \exp\left(-\Delta t\,\hat{u}^\prime\!\!\left(X_{p,\tau+1}\right)\!\right) W^*_{p,\tau+1} - \Delta t \Theta(X_{p,\tau}),
\end{equation}
running in reverse. For each $p$ and $\tau$, this simulation contributes to \eqref{LoBlLe:eq:design} as
\begin{equation}
\frac{\partial B_{p,\tau}}{\partial\hat{u}}^\top\!\!\!(\hat{y}^\prime,\hat{u}^\prime)\hat{\lambda}_{p,\tau} \! = \Delta t W_{p,\tau} \exp\!\left(-\Delta t\,\hat{u}^\prime\!\!\left(\! X_{p,\tau+1}\!\right)\!\right) \!\!
\begin{bmatrix}
	\mathcal{I}_0\!\left(X_{p,\tau+1}\right) &\! \cdots \! & \mathcal{I}_N\!\!\left(X_{p,\tau+1}\right)\!
\end{bmatrix}^\top \! W_{p,\tau}^\ast.
\end{equation}
The total gradient according to \eqref{LoBlLe:eq:design} is then
\begin{equation}
	\frac{\partial\hat{\mathcal{L}}}{\partial\hat{u}}^\top\!\!\!(\hat{y}^\prime, \hat{u}^\prime,\hat{\lambda})  = \nu \Delta x \hat{u}^\prime + \sum_{p=0}^P \sum_{\tau = 1}^T \frac{\partial B_{p,\tau}}{\partial\hat{u}}^\top\!\!\!(\hat{y}^\prime,\hat{u}^\prime)\hat{\lambda}_{p,\tau}.
\end{equation}

\subsection{Reversing particle trajectories}
\label{LoBlLe:sec:adjointchallenges}

We observe that the adjoint simulation \eqref{LoBlLe:eq:heatreversemc} requires the temperature profile $\Theta$, meaning that the full constraint simulation must be completed before starting the adjoint simulation. The adjoint simulation also requires the particle positions $X_{p,\tau}$ in time-reversed order. When $P \times N$ is small, it is possible to simply store the state of each particle at each time step $\tau$, this approach however becomes infeasible as the simulation size increases. Large scale simulations will therefore need to resort to checkpointing, i.e., storing the constraint simulation state at a subset of time steps $\tau$. Such a checkpoint consists of the state vector $\hat{y}_\tau$ and the state of the random number generator at the end of the given time step. During the adjoint simulation, one loads these checkpoints in reverse order and re-computes the required path segments. For information on checkpointing see, e.g.,~\cite[Sec. 5.2]{Naumann2018} and references therein.

While checkpointing avoids an uncontrolled increase in memory usage, it introduces an additional computational cost, which is at least that of simulating the constraint equation, on top of the cost of the adjoint simulation. We aim to avoid the need for checkpointing by re-computing the paths in time reversed order in the adjoint simulation, based on the final state of the constraint simulation. As we can re-write \eqref{LoBlLe:eq:heatrecursion} as an explicit expression for $\hat{y}_\tau$, path generation is reversible if we know the random values involved. To compute these random values we introduce the concept of a reversible random number generator in the next section.

\section{Reversible random number generators}
\label{LoBlLe:sec:rrng}

In this section, we present one of the main contributions of this work, the reversal of the PCG family of pseudorandom number generators~\cite{ONeill2014}. Although we specifically cover PCGs here, we remark that it is typical for the state-transition function of a pseudorandom generator to be a bijection, and thus reversible. This property follows from the fact that a non-bijective state-transition function makes some states more likely than others, thus reducing uniformity. It is however not necessarily the case that the reversal can be done at the same computational cost. Some examples of other candidate reversible generators are (those based on) linear congruential generators~\cite{Brown1994}; counter based generators such as SplitMix~\cite{Steele2014} and Random123~\cite{Salmon2011}; linear feedback shift register based generators such as the Mersenne Twister~\cite{Matsumoto1998} and Xoshiro/Xoroshiro~\cite{Blackman2021}; and chaotic generators such as SFC~\cite{Doty-Humphrey2018}. We first introduce PCGs in Section~\ref{LoBlLe:sec:pcg} and motivate their use as random bit generators. Next in Section~\ref{LoBlLe:sec:reversingpcg}, we reverse PCGs, so that we can request the previous value in a generator's pseudorandom sequence at the same cost as requesting the next value. In Section~\ref{LoBlLe:sec:distributionreversal}, we then describe how we use the reversed generator to sample arbitrary distributions.

\subsection{Permuted congruential generators}
\label{LoBlLe:sec:pcg}

We choose to reverse the family of permuted congruential generators (PCG)~\cite{ONeill2014}. PCGs make internal use of a linear congruential generator (LCG), i.e., generating pseudorandom integer values $\zeta_k$ in the range $[0,m-1]$ through a recurrence
\begin{equation}
	\label{LoBlLe:eq:lcg}
	\zeta_{k+1} = a \zeta_k + c \mod m,
\end{equation}
with $a$ and $c$ integers in the range $[0,m-1]$. However, even for well chosen values $a$, $c$ and $m$,  LCGs have poor statistical properties, as demonstrated by benchmark tests such as TestU01~\cite{LEcuyer2007}. PCGs improve on the poor statistics of the underlying LCG by applying a 1-way parameterized permutation to its state $\zeta_k$ to generate an output value. For technical details and an overview of generator variants, we refer to~\cite{ONeill2014}.

This combination of a simple internal update function and complex output function makes PCG a good candidate for reversal, while still achieving good statistical properties. While the paper~\cite{ONeill2014} was not published in a peer-reviewed venue we have satisfactorily reproduced their TestU01 benchmark results and a number of our runs can be found together with our code referenced in Section~\ref{LoBlLe:sec:experiments}. These generators are also widely used in the scientific computing community, being, e.g., the current default generator in NumPy~\cite{NumPyDevelopers2022}. While they are not suitable as cryptographic generators~\cite{Bouillaguet2020}, and some seed pairs have been found that induce measurable correlation between sequences, see the discussion in~\cite{ONeill2018,Vigna2018}, to our knowledge, no issues have been found in PCGs counterindicating their use in scientific computing.

\subsection{Reversing PCG}
\label{LoBlLe:sec:reversingpcg}

We reverse a PCG by reversing the underlying LCG. That is, we invert the recursion \eqref{LoBlLe:eq:lcg} going from state $\zeta_{k+1}$ to state $\zeta_k$. The 1-way output function can then be applied to these internal state values just as in the forward sequence generation. The problem of reversing the PCG thus reduces to that of reversing the update function \eqref{LoBlLe:eq:lcg}, i.e., performing an update
\begin{equation}
	\label{LoBlLe:eq:reversedlcg}
	\zeta_k = a^{-1} ( \zeta_{k+1} - c ) \mod m,
\end{equation}
with $a^{-1}$ the inverse of $a \mod m$, which is uniquely defined as $a^{-1}\equiv a^{m-2} \mod m$, if the LCG has a period of length $m$. This follows from the Hull-Dobell theorem~\cite[Thm. 1]{Hull1962}, which states that $a$ and $m$ must be co-prime for an LCG to have full period.

One could also compute the previous value by the wraparound of the periodic sequence. However, this computation has a cost equivalent to $\mathcal{O}(\log m)$ generator steps, making reversal very expensive. Instead, we therefore pre-compute the constant $a^{-1} \equiv a^{m-2} \mod m$ and reply the recursion \eqref{LoBlLe:eq:reversedlcg}, which has the same computational cost as \eqref{LoBlLe:eq:lcg}. As previously mentioned, our reversible PCG implementaton requires access to the internal state. As such, our \verb|C++| implementation directly inherits from the existing PCG class. Our implementation can be found both together with the code for the experiments in Section~\ref{LoBlLe:sec:experiments}, and in a pull request to the PCG repository\footnote{\url{https://github.com/imneme/pcg-cpp/pull/77}}.

\subsection{Sampling distributions in reverse}
\label{LoBlLe:sec:distributionreversal}

With our reversible PCG, producing sequences of uniformly distributed 64 bit integers in both forward and reverse mode, we want to generate IEEE-754 double precision values from various distributions. We first consider the uniform distribution $\mathcal{U}([0,1))$. Here, it is important to note that it is not possible to define a bijective mapping between 64 bit integers and floating point values that maintains uniformity. This follows from the non-uniformity of floating-point numbers. Here, we take a pragmatic approach, in which we restrict the integers to 53 bits of information, matching the number of floating-point mantissa-bits (52 physical bits and one implicitly defined), thus discarding 11 bits. Discarding these bits ensures that no rounding occurs in the type conversion, thus maintaining uniformity at the cost of additional quantization. In \verb|C++| this transformation is given by
\begin{lstlisting}
static inline double float64(std::uint64_t x)
{
	return (x >> 11) * 0x1.0p-53;
}
\end{lstlisting}
If the resulting quantization would be unacceptable, alternative methods exist that can generate all 64-bit floating point values by using more input bits, see e.g.~\cite{Campbell2014}.

Once we can sample the uniform distribution, there are two common approaches to sampling other distributions. These are applying analytical transformations through the inverse cumulative density function and using accept-reject sampling. We now discuss reverse sampling in both of these cases. We use the exponential and normal distributions as examples, while remarking that at least one of these approaches can be applied to any commonly used distribution in scientific computing.

\subsubsection{Inverse CDF transformations}

For some distributions, a simple analytical expression exists for converting uniform values to samples from the desired distribution. The exponential distribution $\mathcal{E}(R)$, describing the time between events in a Poisson point process with event rate $R$, is one such example. Given the cumulative density function for some $\eta \sim \mathcal{E}(R)$, we set 
\begin{equation}
	\phi = \text{CDF}_{\mathcal{E}}(\eta) = 1-e^{-R \eta},
\end{equation}
and observe that $\phi \sim \mathcal{U}([0,1))$. We then know that  the inverse function of the CDF
\begin{equation}
	\eta = \text{CDF}^{-1}_{\mathcal{E}} (\phi) = -\frac{\ln (1-\phi)}{R}
\end{equation}
transforms a value $\phi\sim\mathcal{U}([0,1))$ to a value $\eta \sim \mathcal{E}(R)$. This approach can be applied to any distribution for which a simple expression for the inverse CDF exists and does not rely on the order in which a sequence of values $\phi$ is generated. It however risks introducing quantization effects in low-probability regions~\cite[Ch.~4]{Gentle2003}.

\subsubsection{Accept-reject sampling}

When the inverse CDF is too complex for use as a transformation or quantization in, e.g., the tails, is unacceptable, an alternative is to use an accept-reject algorithm. This mean that we sample pseudorandom values from a simpler distribution and selectively discard samples so that the non-discarded samples have the desired distribution. Here, we consider the case of sampling the standard normal distribution $\mathcal{N}(0,1)$ using the Ziggurat algorithm~\cite{Marsaglia2000}.

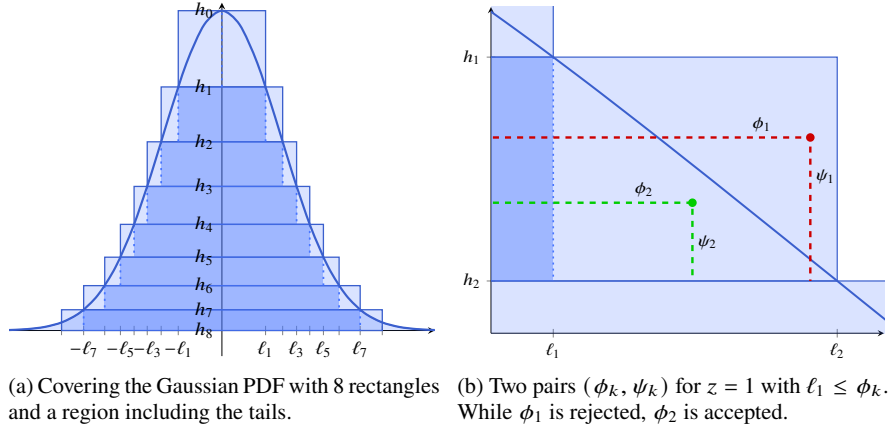
\begin{figure}
	\centering
		\begin{minipage}{\textwidth}
		\subfloat[Covering the Gaussian PDF with 8 rectangles and a region including the tails. \label{LoBlLe:fig:zigguratfull}]{\resizebox{0.49\textwidth}{!}{\input{ziggurat}}}
		\hfill
		\subfloat[Two pairs $(\phi_k,\psi_k)$ for $z=1$ with $\ell_1 \leq \phi_k$. While $\phi_1$ is rejected, $\phi_2$ is accepted. \label{LoBlLe:fig:zigguratzoom}]{\resizebox{0.49\textwidth}{!}{\input{zigguratzoom}}}
	\end{minipage}
\caption{A sketch of the Ziggurat accept-reject algorithm. Each horizontal light-blue region has the same area, while dark-blue regions mark the range of values for $\phi_k$ which are always accepted.} \label{LoBlLe:fig:ziggurat}
\end{figure}

The Ziggurat algorithm approximates the Bell curve of the normal PDF by a set of $Z$ rectangles with horizontal span $[-\ell_{z+1}, \ell_{z+1}]$ and vertical span $(h_{z+1}, h_z]$, with $z=0,\dots, Z-1$ and $Z$ a power of 2. These rectangles cover the majority of the curve as sketched in Figure~\ref{LoBlLe:fig:zigguratfull} and each have the same area as the intersection of the Bell curve and the band $[0,h_{Z-1}]$. Given a value $\zeta_k$ with 64 random bits, the algorithm uses $\log_2(Z)$ of the bits to select an index $z_k \in [0,Z-1]$. The remaining bits of $\zeta_k$ are used to generate a random value $\phi_k$, according to one of two cases:
\begin{itemize}
	\item If $z_k \neq Z-1$, then $\phi_k\sim \mathcal{U}([-\ell_{z_k+1},\ell_{z_k+1}])$. If $|\phi_k| < \ell_{z_k}$, $\phi_k$ is accepted. Otherwise, we sample an additional value $\psi_k \sim \mathcal{U}((h_{z_k+1}, h_{z_k}])$. If $\psi_k < \exp(-\frac{\phi_k^2}{2})$ then $\phi_k$ is accepted, otherwise $\phi_k$ is discarded and the algorithm repeats for a new value $\zeta_{k+1}$. We sketch this step in Figure~\ref{LoBlLe:fig:zigguratzoom}.
	\item If $z_k=Z-1$, $\phi_k\sim \mathcal{U}([-\ell_{Z},\ell_{Z}])$. If $|\phi_k| < \ell_{Z-1}$, then $\phi_k$ is returned. Otherwise, we sample the tails beyond the final rectangle with the Marsaglia method~\cite{Marsaglia1964a}, which repeatedly converts pairs of uniform values into a proposed normal value from the tails to be accepted or rejected.
\end{itemize}
The Ziggurat algorithm, follows the \say{make the common case fast} principle. For example, an implementation with $Z=128$ accepts the initially proposed $\phi_k$ in 98.78\% of cases and performs even better as the number of rectangles increases~\cite{Marsaglia2000}.

Depending on the value of $\zeta_k$, the Ziggurat algorithm requires either 0, 1 or $1+2i$, $i\in\mathbb{N}$, uniform values before sampling a next value $\zeta_{k+1}$. If all uniform values are drawn form the same sequence, then it is not possible to determine which elements of the sequence correspond with values $\zeta_k$. We therefore take the approach of only generating the values $\zeta_k$ from our reversible PCG implementation. If more uniform values are needed to generate a proposal, we then use $\zeta_k$ to seed a second pseudorandom number generator from which we get these additional values.

In the most common case, i.e., $z_k \neq Z-1$ and $|\phi_k| <\ell_{z_k}$, no other random values are needed and the performance of our reversible implementation remains unchanged. In other cases, this seeding induces a small amount of extra overhead on top of the Ziggurat algorithm. We deem this small overhead acceptable in the small percentage of cases where such additional values are necessary. To avoid inducing correlations between these additional uniform values and future values of the $\zeta$ sequence, we refrain from using a PCG to generate these additional values. Instead we use the \texttt{xoshiro256+} generator~\cite{Blackman2021} for this purpose.

This second generator \say{trick} makes our accept-reject strategy more general than that proposed in~\cite{Yoginath2018}, but introduces a potentially biased values. This bias is caused by repeated seeding of the second generator, with fewer possible seed values than can be represented in the generator state. Following the design of the \verb|C++| standard library, we separate bit generation from distribution sampling, meaning that we cannot access the PCG's internal state but only the 64 output bits. Inside \texttt{xoshiro256+}, SplitMix~\cite{Steele2014} is used to populate its 256 bit state from a 64 bit seed. Additionally, we only use seed values where $z_k = Z-1$ or $|\phi_k| \geq \ell_{z_k}$. As a consequence, the second generator cannot reach all possible states. Despite these considerations, the values produced don't consistently fail any tests in TestU01 BigCrush after conversion to uniform values using a CDF transformation. A number of our runs can be found together with our code referenced in Section~\ref{LoBlLe:sec:experiments}.

Though we only consider the Ziggurat algorithm here, our seeding approach works for any accept-reject algorithm with uncorrelated proposals. Rather than seeding a second generator, one an equally general bias-free approach would be to use a splittable generator such as SplitMix~\cite{Steele2014} or using multiple streams of same generator. These approaches come at the cost of software modularity. One can also reverse an arbitrary accept-reject algorithm without bias and without breaking modularity by using multiple reversible generators. At most, one needs a generator for each random variable occurring in a branch-condition affecting the number of random variables needed to generate a proposal. The Ziggurat algorithm can, e.g., be reversed with two reversible generators, one providing values $\zeta_k$ and another providing the values $\psi_k$ and the pairs used in the Marsaglia algorithm, as the required access pattern of the values from the second sequence is known once $\zeta_k$ is determined.

\section{Simulation results}
\label{LoBlLe:sec:experiments}

We now demonstrate the efficacy of our reversible generator approach. First, in Section~\ref{LoBlLe:sec:generatorresults}, we consider the performance of the reversible generator. Next, in Section~\ref{LoBlLe:sec:cyclicresults} we consider the problem from Section~\ref{LoBlLe:sec:1Dheat}. The code for producing our results is available at \url{gitlab.kuleuven.be/numa/public/mcqmc\_2022\_rrng}, together with a Docker configuration with the software dependencies. All timings were produced in a Podman container on a 2021 M1 MacBook Pro with 16GB of RAM.

\subsection{Generator efficiency}
\label{LoBlLe:sec:generatorresults}

\begin{figure}
	\centering
	
	\subfloat[Uniform distribution\label{LoBlLe:fig:reversaluniform}]{
		\hspace{-14pt}
		\resizebox{0.49\textwidth}{!}{
			\begin{tikzpicture}
				\begin{axis}[
					xlabel={Random values generated},
					ylabel={Runtime $[s]$},
					xmode=log,
					xmax=1.5e8,
					xmin=7e2,
					ymax=1,
					ymin=1e-6,
					ymode=log,
					legend pos=north west
					]
					\addplot[thick, RoyalBlue2, mark=none] table [x=samples, y=uniformForward, col sep=comma] {"figure2.csv"} ;
					\addplot[very thick, dotted, Red4, mark=none, mark options={solid}] table [x=samples, y=uniformBackward, col sep=comma] {"figure2.csv"} ;
					\legend{Forward, Reverse}
				\end{axis}
		\end{tikzpicture}}
		\hspace{14pt}}

	\subfloat[Normal distribution\label{LoBlLe:fig:reversalnormal}]{
		\hspace{-14pt}
		\resizebox{0.49\textwidth}{!}{
			\begin{tikzpicture}
				\begin{axis}[
					xlabel={Random values generated},
					ylabel={Runtime $[s]$},
					xmode=log,
					xmax=1.5e8,
					xmin=7e2,
					ymax=1,
					ymin=1e-6,
					ymode=log,
					legend pos=north west
					]
					\addplot[thick, RoyalBlue2, mark=none] table [x=samples, y=normalForward, col sep=comma] {"figure2.csv"} ;
					\addplot[very thick, dotted, Red4, mark=none, mark options={solid}] table [x=samples, y=normalBackward, col sep=comma] {"figure2.csv"} ;
					\legend{Forward, Reverse}
				\end{axis}
		\end{tikzpicture}}
		\hspace{14pt}}
	\subfloat[Exponential distribution\label{LoBlLe:fig:reversalexponential}]{
		\hspace{-14pt}
		\resizebox{0.49\textwidth}{!}{
			\begin{tikzpicture}
				\begin{axis}[
					xlabel={Random values generated},
					ylabel={Runtime $[s]$},
					xmode=log,
					xmax=1.5e8,
					xmin=7e2,
					ymax=1,
					ymin=1e-6,
					ymode=log,
					legend pos=north west
					]
					\addplot[thick, RoyalBlue2, mark=none] table [x=samples, y=exponentialForward, col sep=comma] {"figure2.csv"} ;
					\addplot[very thick, dotted, Red4, mark=none, mark options={solid}] table [x=samples, y=exponentialBackward, col sep=comma] {"figure2.csv"} ;
					\legend{Forward, Reverse}
				\end{axis}
		\end{tikzpicture}}
		\hspace{14pt}}
	
	\caption{Comparing forward and reverse mode timings for different random distributions. We see that both modes attain the same performance for all three considered distributions.}\label{LoBlLe:fig:reversal}
\end{figure}
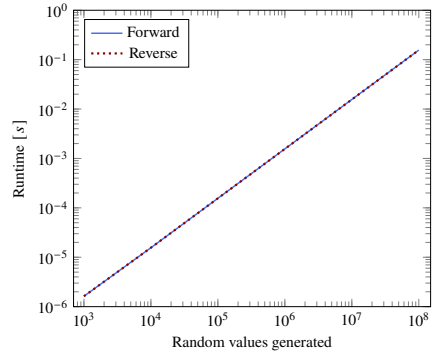
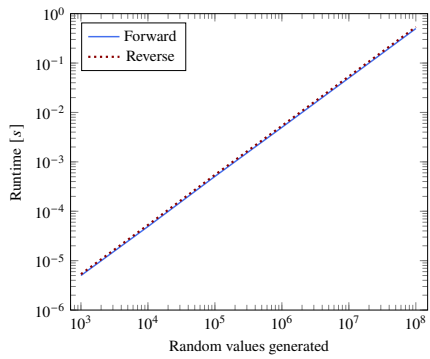
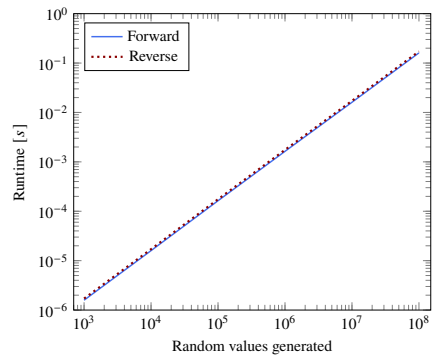

To test the efficiency of our reversed generator we time how long our \verb|C++| implementation needs to generate $K$ random values in both forward and reverse mode for different seeds, with $K$ varying between $10^3$ and $10^8$. To produce reliable timings, we perform 55 timings for each $K$, discarding the first 5 timings and reporting the minimum of the remaining 50. The resulting timings for \verb|-O3| are shown in Figure~\ref{LoBlLe:fig:reversal}.

In Figure~\ref{LoBlLe:fig:reversaluniform}, we see that our reverse-mode extension of PCG attains the same performance as the (unaltered) forward implementation. We have thus reversed the uniform generator, at the same computational cost as the forward generation. In addition, Figures~\ref{LoBlLe:fig:reversalnormal} and~\ref{LoBlLe:fig:reversalexponential} show that both forward and reverse modes for our non-uniform generators attain the same performance. While the cost of sampling the exponential distribution is practically the same as sampling the uniform distribution, sampling the normal distribution values is approximately twice as expensive.

\subsection{Cooling a 1D heat equation}
\label{LoBlLe:sec:cyclicresults}

\begin{figure}
	\centering
	\resizebox{0.75\textwidth}{!}{
		\begin{tikzpicture}
		\begin{axis}[
			xlabel={$x$},
			ylabel={$\hat{\theta}(x,t)$, $\hat{u}(x)$},
			xmin=0,
			xmax=10,
			ymin=-5,
			ymax=165,
			legend style={at={(1.05,0.5)},anchor=west}
			]
			\addplot[thick, RoyalBlue2, mark=none] table [x index=1, y index=2] {\HeatOptimizationTransposed} ;
			\addplot[thick, dashed, Red4, mark=none, mark options={solid}] table [x index=1, y index=3] {\HeatOptimizationTransposed} ;
			\addplot[thick, dotted, DarkOliveGreen4, mark=none, mark options={solid}] table [x index=1, y index=4] {\HeatOptimizationTransposed} ;
			
			\legend{{$\hat{\theta}(x,0)$}, {$\hat{\theta}(x,T)$}, {$\hat{u}(x)$}}
		\end{axis}
	\end{tikzpicture}
	}
	\caption{The temperature profile and optimal cooling rate for the 1D heat equation problem.\label{LoBlLe:fig:cyclicheatsolution}}
\end{figure}
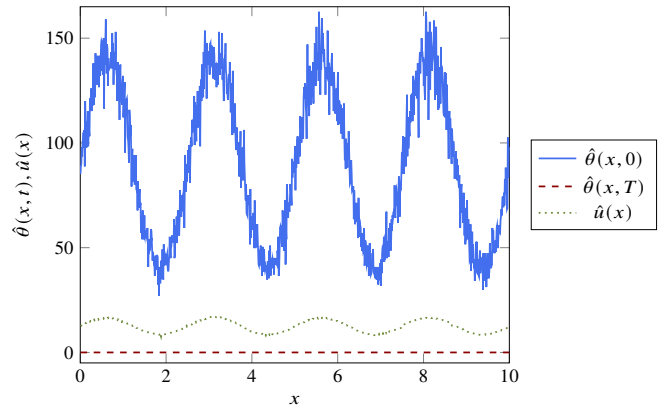

We now consider the simulating the 1D heat equation problem from Section~\ref{LoBlLe:sec:1Dheat}. The periodic domain allows us to focus on reversing the path dynamics themselves, without the additional complexities induced by more complex boundary conditions. As an initial condition, we set $y(x,0) = 50 \sin (8 \pi x/L) + 90$ and set the regularization parameter $\nu=1$, the domain length $L=10$ and end time $T=1$. We set the spatial grid size to $\Delta x = 0.01$ and the time step size to $\Delta t = 0.001$.  We show the initial and final conditions, as well as the optimal computed $\hat{u}(x)$ after convergence in Figure~\ref{LoBlLe:fig:cyclicheatsolution}.

\pgfplotstableread[col sep=comma]{"figure3.csv"}{\HeatOptimization}
\pgfplotstabletranspose{\HeatOptimizationTransposed}{\HeatOptimization}

\pgfplotstableread[col sep=comma]{"figure4a.csv"}{\HeatTimings}
\pgfplotstabletranspose{\HeatTimingsTransposed}{\HeatTimings}
\pgfplotstableread[col sep=comma]{"figure4b.csv"}{\HeatMemory}
\pgfplotstabletranspose{\HeatMemoryTransposed}{\HeatMemory}
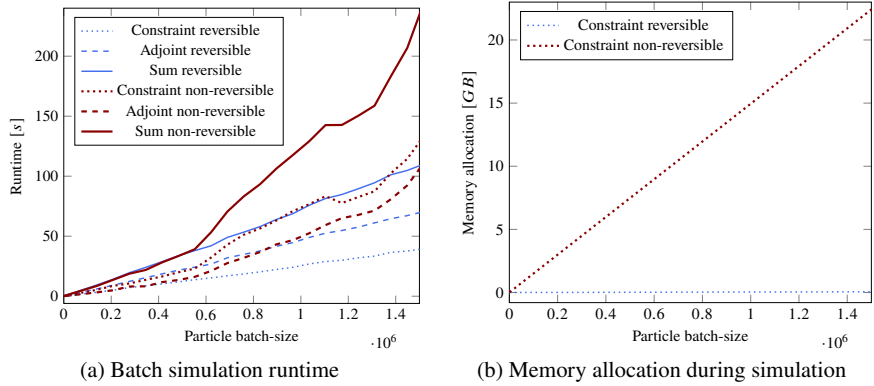
\begin{figure}
	\subfloat[Batch simulation runtime\label{LoBlLe:fig:cyclicheatanalysisruntime}]{
		\resizebox{0.48\textwidth}{!}{
			\begin{tikzpicture}
				\begin{axis}[
					xlabel={Particle batch-size},
					ylabel={Runtime $[s]$},
					xmin=0,
					xmax=1.5e6,
					ymin=-6,
					ymax=240,
					legend pos=north west
					]
					\addplot[thick, dotted, RoyalBlue2, mark=none] table [x index=1, y index=2] {\HeatTimingsTransposed} ;
					\addplot[thick, dashed, RoyalBlue2, mark=none, mark options={solid}] table [x index=1, y index=3] {\HeatTimingsTransposed} ;
					\addplot[thick, RoyalBlue2, mark=none, mark options={solid}] table [x index=1, y index=4] {\HeatTimingsTransposed} ;
					\addplot[very thick, dotted, Red4, mark=none, mark options={solid}] table [x index=1, y index=5] {\HeatTimingsTransposed} ;
					\addplot[very thick, dashed, Red4, mark=none, mark options={solid}] table [x index=1, y index=6] {\HeatTimingsTransposed} ;
					\addplot[very thick, Red4, mark=none, mark options={solid}] table [x index=1, y index=7] {\HeatTimingsTransposed} ;
					
					\legend{Constraint reversible, Adjoint reversible, Sum reversible, Constraint non-reversible, Adjoint non-reversible, Sum non-reversible}
				\end{axis}
			\end{tikzpicture}
	}}
	\hfill
	\subfloat[Memory allocation during simulation]{
		\resizebox{0.48\textwidth}{!}{
			\begin{tikzpicture}
				\begin{axis}[
					xlabel={Particle batch-size},
					ylabel={Memory allocation $[GB]$},
					xmin=0,
					xmax=1.5e6,
					ymin=-0.8,
					ymax=23,
					y filter/.code={\pgfmathparse{#1/1024/1024/1024}\pgfmathresult},
					legend pos=north west
					]
					\addplot[thick, dotted, RoyalBlue2, mark=none] table [x index=1, y index=2] {\HeatMemoryTransposed} ;
					\addplot[very thick, dotted, Red4, mark=none, mark options={solid}] table [x index=1, y index=5] {\HeatMemoryTransposed} ;
					\legend{Constraint reversible, Constraint non-reversible}
				\end{axis}
			\end{tikzpicture}
	}}
	\caption{The increase in runtime and memory usage of a reversible and non-reversible constraint and corresponding adjoint simulations as the number of particles increases. As memory requirements approach the available RAM on the system (16GB), the non-reversible simulation's runtime increases drastically, while the reversible implementation maintains its linear scaling.\label{LoBlLe:fig:cyclicheatanalysis}}
\end{figure}

We now compare a classic implementation of this problem, storing paths in memory, to one where we re-compute the paths in the adjoint simulation using a reversible generator. We compare both approaches for varying batch-sizes, i.e., the number of particles simulated simultaneously, in terms of computational cost and memory usage. In Figure~\ref{LoBlLe:fig:cyclicheatanalysisruntime}, we see that both implementations perform comparably for smaller batch-sizes. However, once the memory requirements of the algorithm approaches the available 16GB of RAM on the computer, the runtime of the classic implementation increases sharply, while the cost of the reversible implementation continues to grow at an unchanged rate. At this point in the plot, the classic implementation starts writing to disk, causing severe slowdowns in the computation.

While one can avoid this catastrophic blow-up in computational cost through check-pointing, we avoid implementing and fine-tuning a check-pointing code here, as we know that it by definition has a computational cost higher than that of the reversible implementation due to the required recomputation. We however assume that a lower bound for the cost of a checkpointed implementation is given by two reversible constraint simulations and a non-reversible adjoint simulation.

\section{Conclusions}
\label{LoBlLe:sec:conclusions}

We presented a reversible extension to PCG pseudorandom number generators. We also presented strategies for reversing arbitrary distributions using either transformation and accept-reject approaches. We demonstrated that this reversible approach allows us to generate a given sequence of values with the same computational cost both forward from the initial generator state and backwards from the final state.

We then applied this reversible generation to the setting of discrete adjoint-based PDE-constrained optimization with Monte Carlo. We demonstrated on a test problem involving the 1D heat equation that our reversible approach avoids the high memory costs of storing the paths of the constraint equation, while still allowing for a simple implementation of the adjoint solver. We therefore present this approach as an alternative to more complex and more expensive ad-hoc checkpointing approaches.

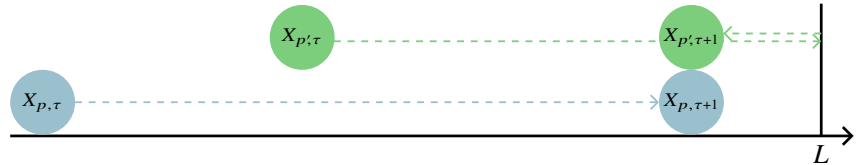
\begin{figure}
	\centering
	\resizebox{\textwidth}{!}{
	\input{reflections}
	}
	\caption{Demonstration of the challenge of reflective boundaries. Both particles $X_p$ and $X_{p^\prime}$ have the same position at time $\tau+1$ after moving the same distance to the right, even though they have different positions at time $\tau$.\label{LoBlLe:fig:reflections}}
\end{figure}

Although path-reversal is conceptually straightforward, we foresee certain cases were non-injective operations may present challenges in more complex simulations. One such example is reflective boundaries, where, as demonstrated by Figure~\ref{LoBlLe:fig:reflections}, there are two locations from which a particle can reach a given end-position given a known positional increment. Such issues are often down to implementation, rather than fundamental. The given example can, e.g., be avoided by the introduction of a reflection vector as an additional particle state, which is updated on each reflection. In 1D, this reflection vector based solution, can be interpreted as a simulation on a domain of length $2L$ with periodic boundaries.

In future work, we plan to apply this reversible approach to more complex and higher-dimensional PDEs, where Monte Carlo methods are highly attractive. An example is fusion reactor design with kinetic equations as a constraint~\cite{Dekeyser2018}.

\begin{acknowledgement}
This paper is based on the master's theses~\cite{Blondeel2022, Lee2022, Vanroye2019}. We thank the anonymous reviewer for their detailed comments on our discussion of random number generation and Mike Giles for pointing out a mathematical error in the first version of this paper. We thank Michael Mascagni for his input on suitability of different random number generators for this work. We also thank Ignace Bossuyt, Vince Maes and Zhirui Tang for providing feedback on our experiments' reproducibility. During this work, Emil Løvbak and Andreas Van Barel were funded by the Research Foundation - Flanders (FWO) under, respectively, fellowship numbers 1SB1919N/1SB1921N and 11E1518N. This work has been carried out within the framework of the EUROfusion Consortium, funded by the European Union via the Euratom Research and Training Programme (Grant Agreement No 101052200 — EUROfusion). Views and opinions expressed are however those of the author(s) only and do not necessarily reflect those of the European Union or the European Commission. Neither the European Union nor the European Commission can be held responsible for them.
\end{acknowledgement}

\bibliographystyle{spmpsci}
\bibliography{references}
\end{document}

%% file: ziggurat.tex
\begin{tikzpicture} 
	\begin{axis} [
		axis lines = middle,
		ytick style={draw=none},
		ymax=1.025,
		ymin =-0.08,
		yticklabels={,,},
		xtick={2.712053022, 2.3383716982418785,1.981904936397471, 1.7165081257747172, 1.4853586756417962, 1.2629701985297244, 1.0273863717794265, 0.738368917975931,  -2.712053022, -2.3383716982418785,-1.981904936397471, -1.7165081257747172, -1.4853586756417962, -1.2629701985297244, -1.0273863717794265, -0.738368917975931},
		xticklabels={,$\ell_7$,,$\ell_5$,,$\ell_3$,,$\ell_1$,,$-\ell_7$,,$-\ell_5$,,$-\ell_3$,,$-\ell_1$}
		]
		\draw[semithick, draw=RoyalBlue3, fill=RoyalBlue1,fill opacity=0.2] (-0.738368917975931 ,1) rectangle (0.738368917975931 ,{exp{-0.738368917975931 *0.738368917975931 /2}});
		\draw[semithick, draw=RoyalBlue3, fill=RoyalBlue1,fill opacity=0.2] (-1.0273863717794265,{exp{-0.738368917975931 *0.738368917975931 /2}}) rectangle (1.0273863717794265,{exp{-1.0273863717794265*1.0273863717794265/2}});
		\draw[semithick, draw=RoyalBlue3,  fill=RoyalBlue1,fill opacity=0.2] (-1.2629701985297244,{exp{-1.0273863717794265*1.0273863717794265/2}}) rectangle (1.2629701985297244,{exp{-1.2629701985297244*1.2629701985297244/2}});
		\draw[semithick, draw=RoyalBlue3,  fill=RoyalBlue1,fill opacity=0.2] (-1.4853586756417962,{exp{-1.2629701985297244*1.2629701985297244/2}}) rectangle (1.4853586756417962,{exp{-1.4853586756417962*1.4853586756417962/2}});
		\draw[semithick, draw=RoyalBlue3,  fill=RoyalBlue1,fill opacity=0.2] (-1.7165081257747172,{exp{-1.4853586756417962*1.4853586756417962/2}}) rectangle (1.7165081257747172,{exp{-1.7165081257747172*1.7165081257747172/2}});
		\draw[semithick, draw=RoyalBlue3,  fill=RoyalBlue1,fill opacity=0.2] (-1.981904936397471,{exp{-1.7165081257747172*1.7165081257747172/2}}) rectangle (1.981904936397471,{exp{-1.981904936397471*1.981904936397471/2}});
		\draw[semithick, draw=RoyalBlue3,  fill=RoyalBlue1,fill opacity=0.2] (-2.3383716982418785,{exp{-1.981904936397471*1.981904936397471/2}}) rectangle (2.3383716982418785,{exp{-2.3383716982418785*2.3383716982418785/2}});
		\draw[semithick, draw=RoyalBlue3,  fill=RoyalBlue1,fill opacity=0.2] (-2.712053022,{exp{-2.3383716982418785*2.3383716982418785/2}}) rectangle (2.712053022,0);
		
		\draw[thick, dash pattern=on 1pt off 3pt, draw=RoyalBlue4, draw opacity=0, fill=RoyalBlue1,fill opacity=0.4] (-0.738368917975931 ,{exp{-0.738368917975931 *0.738368917975931 /2}}) rectangle (0.738368917975931 ,{exp{-1.0273863717794265 *1.0273863717794265 /2}});
		\draw[thick, dash pattern=on 1pt off 3pt, draw=RoyalBlue4, draw opacity=0, fill=RoyalBlue1,fill opacity=0.4] (-1.0273863717794265 ,{exp{-1.0273863717794265 *1.0273863717794265 /2}}) rectangle (1.0273863717794265 ,{exp{-1.2629701985297244 *1.2629701985297244 /2}});
		\draw[thick, dash pattern=on 1pt off 3pt, draw=RoyalBlue4, draw opacity=0, fill=RoyalBlue1,fill opacity=0.4] (-1.2629701985297244 ,{exp{-1.2629701985297244 *1.2629701985297244 /2}}) rectangle (1.2629701985297244 ,{exp{-1.4853586756417962 *1.4853586756417962 /2}});
		\draw[thick, dash pattern=on 1pt off 3pt, draw=RoyalBlue4, draw opacity=0, fill=RoyalBlue1,fill opacity=0.4] (-1.4853586756417962 ,{exp{-1.4853586756417962 *1.4853586756417962 /2}}) rectangle (1.4853586756417962 ,{exp{-1.7165081257747172 *1.7165081257747172 /2}});
		\draw[thick, dash pattern=on 1pt off 3pt, draw=RoyalBlue4, draw opacity=0, fill=RoyalBlue1,fill opacity=0.4] (-1.7165081257747172 ,{exp{-1.7165081257747172 *1.7165081257747172 /2}}) rectangle (1.7165081257747172 ,{exp{-1.9819049363974715 *1.981904936397471 /2}});
		\draw[thick, dash pattern=on 1pt off 3pt, draw=RoyalBlue4, draw opacity=0, fill=RoyalBlue1,fill opacity=0.4] (- 1.981904936397471  ,{exp{- 1.981904936397471  * 1.981904936397471  /2}}) rectangle (1.981904936397471  ,{exp{-2.3383716982418785 *2.3383716982418785 /2}});
		\draw[thick, dash pattern=on 1pt off 3pt, draw=RoyalBlue4, draw opacity=0, fill=RoyalBlue1,fill opacity=0.5] (- 2.3383716982418785  ,{exp{- 2.3383716982418785 * 2.3383716982418785  /2}}) rectangle (2.3383716982418785 , 0);
		
		\draw[thick, dash pattern=on 1pt off 2pt, draw=RoyalBlue1] (2.3383716982418785  ,{exp{- 2.3383716982418785 * 2.3383716982418785  /2}}) -- (2.3383716982418785 , 0);
		\draw[thick, dash pattern=on 1pt off 2pt, draw=RoyalBlue1] (-2.3383716982418785  ,{exp{- 2.3383716982418785 * 2.3383716982418785  /2}}) -- (-2.3383716982418785 , 0);
		\draw[thick, dash pattern=on 1pt off 2pt, draw=RoyalBlue1] (1.981904936397471  ,{exp{- 1.981904936397471 * 1.981904936397471  /2}}) -- (1.981904936397471 , {exp{- 2.3383716982418785 * 2.3383716982418785  /2}});
		\draw[thick, dash pattern=on 1pt off 2pt, draw=RoyalBlue1] (-1.981904936397471  ,{exp{- 1.981904936397471 * 1.981904936397471  /2}}) -- (-1.981904936397471 , {exp{- 2.3383716982418785 * 2.3383716982418785  /2}});
		\draw[thick, dash pattern=on 1pt off 2pt, draw=RoyalBlue1] (1.7165081257747172  ,{exp{- 1.7165081257747172 * 1.7165081257747172  /2}}) -- (1.7165081257747172 , {exp{- 1.981904936397471 * 1.981904936397471  /2}});
		\draw[thick, dash pattern=on 1pt off 2pt, draw=RoyalBlue1] (-1.7165081257747172  ,{exp{- 1.7165081257747172 * 1.7165081257747172  /2}}) -- (-1.7165081257747172 , {exp{- 1.981904936397471 * 1.981904936397471  /2}});
		\draw[thick, dash pattern=on 1pt off 2pt, draw=RoyalBlue1] (1.4853586756417962  ,{exp{- 1.4853586756417962 * 1.4853586756417962  /2}}) -- (1.4853586756417962 , {exp{- 1.7165081257747172 * 1.7165081257747172  /2}});
		\draw[thick, dash pattern=on 1pt off 2pt, draw=RoyalBlue1] (-1.4853586756417962  ,{exp{- 1.4853586756417962 * 1.4853586756417962  /2}}) -- (-1.4853586756417962 , {exp{- 1.7165081257747172 * 1.7165081257747172  /2}});
		\draw[thick, dash pattern=on 1pt off 2pt, draw=RoyalBlue1] (1.2629701985297244  ,{exp{- 1.2629701985297244 * 1.2629701985297244  /2}}) -- (1.2629701985297244 , {exp{- 1.4853586756417962 * 1.4853586756417962  /2}});
		\draw[thick, dash pattern=on 1pt off 2pt, draw=RoyalBlue1] (-1.2629701985297244  ,{exp{- 1.2629701985297244 * 1.2629701985297244  /2}}) -- (-1.2629701985297244 , {exp{- 1.4853586756417962 * 1.4853586756417962  /2}});
		\draw[thick, dash pattern=on 1pt off 2pt, draw=RoyalBlue1] (1.0273863717794265  ,{exp{- 1.0273863717794265 * 1.0273863717794265  /2}}) -- (1.0273863717794265 , {exp{- 1.2629701985297244 * 1.2629701985297244  /2}});
		\draw[thick, dash pattern=on 1pt off 2pt, draw=RoyalBlue1] (-1.0273863717794265  ,{exp{- 1.0273863717794265 * 1.0273863717794265  /2}}) -- (-1.0273863717794265 , {exp{- 1.2629701985297244 * 1.2629701985297244  /2}});
		\draw[thick, dash pattern=on 1pt off 2pt, draw=RoyalBlue1] (0.738368917975931  ,{exp{- 0.738368917975931 * 0.738368917975931  /2}}) -- (0.738368917975931 , {exp{- 1.0273863717794265 * 1.0273863717794265  /2}});
		\draw[thick, dash pattern=on 1pt off 2pt, draw=RoyalBlue1] (-0.738368917975931  ,{exp{- 0.738368917975931 * 0.738368917975931  /2}}) -- (-0.738368917975931 , {exp{- 1.0273863717794265 * 1.0273863717794265  /2}});
		\draw[thick, dash pattern=on 1pt off 2pt, draw=RoyalBlue1] (0 , {exp{- 0.738368917975931 * 0.738368917975931  /2}}) -- (0,1);		
		
		\addplot[name path=fright,domain=2.3383716982418785:3.6,blue] {exp{-x*x/2}};
		\path[name path=axisright] (axis cs:2.3383716982418785,0) -- (axis cs:3.6,0);
		\addplot [
		thick,
		color=blue,
		fill=RoyalBlue1, 
		fill opacity=0.2
		]
		fill between[
		of=fright and axisright,
		soft clip={domain=2.3383716982418785:3.6},
		];
		
		\addplot[name path=fright,domain=-3.6:-2.3383716982418785,blue] {exp{-x*x/2}};
		\path[name path=axisright] (axis cs:-3.6,0) -- (axis cs:-2.3383716982418785,0);
		\addplot [
		thick,
		color=blue,
		fill=RoyalBlue1, 
		fill opacity=0.2
		]
		fill between[
		of=fright and axisright,
		soft clip={domain=-3.6:-2.3383716982418785},
		];
		
		\addplot [RoyalBlue3, line width = 1, smooth, domain=-3.6:3.6] {exp{-x*x/2}};
		
		\draw (-0.3,1) node {$h_0$};
		\draw (-0.3,{exp{- 0.738368917975931 * 0.738368917975931  /2}}) node {$h_1$};
		\draw (-0.3,{exp{- 1.027386371779426 * 1.027386371779426  /2}}) node {$h_2$};
		\draw (-0.3,{exp{- 1.2629701985297244 * 1.2629701985297244  /2}}) node {$h_3$};
		\draw (-0.3,{exp{- 1.4853586756417962 * 1.4853586756417962  /2}}) node {$h_4$};
		\draw (-0.3,{exp{- 1.7165081257747172 * 1.7165081257747172 /2}}) node {$h_5$};
		\draw (-0.3,{exp{- 1.981904936397471 * 1.981904936397471 /2}}) node {$h_6$};
		\draw (-0.3,{exp{- 2.3383716982418785 * 2.3383716982418785 /2}}) node {$h_7$};
		\draw (-0.3,0) node {$h_8$};
		
		\end{axis} 
\end{tikzpicture}


%% file: zigguratzoom.tex
\begin{tikzpicture} 
	\begin{axis} [
		axis lines = middle,
		ymax=0.8,
		ymin =0.55,
		yticklabels={,,},
		xtick={1.0273863717794265, 0.738368917975931},
		xticklabels={$\ell_2$,$\ell_1$},
		ytick={{exp{- 1.027386371779426 * 1.027386371779426  /2}},{exp{- 0.738368917975931 * 0.738368917975931  /2}}},
		yticklabels={$h_2$,$h_1$}
		]
		\draw[semithick, draw=RoyalBlue3, fill=RoyalBlue1,fill opacity=0.2] (-0.738368917975931 ,1) rectangle (0.738368917975931 ,{exp{-0.738368917975931 *0.738368917975931 /2}});
		\draw[semithick, draw=RoyalBlue3, fill=RoyalBlue1,fill opacity=0.2] 	(-1.0273863717794265,{exp{-0.738368917975931 *0.738368917975931 /2}}) rectangle (1.0273863717794265,{exp{-1.0273863717794265*1.0273863717794265/2}});
		\draw[semithick, draw=RoyalBlue3,  fill=RoyalBlue1,fill opacity=0.2] (-1.2629701985297244,{exp{-1.0273863717794265*1.0273863717794265/2}}) rectangle (1.2629701985297244,{exp{-1.2629701985297244*1.2629701985297244/2}});
		
		\draw[thick, dash pattern=on 1pt off 3pt, draw=RoyalBlue4, draw opacity=0, fill=RoyalBlue1,fill opacity=0.4] (-0.738368917975931 ,{exp{-0.738368917975931 *0.738368917975931 /2}}) rectangle (0.738368917975931 ,{exp{-1.0273863717794265 *1.0273863717794265 /2}});
		
		\draw[thick, dash pattern=on 1pt off 2pt, draw=RoyalBlue1] (0.738368917975931  ,{exp{- 0.738368917975931 * 0.738368917975931  /2}}) -- (0.738368917975931 , {exp{- 1.0273863717794265 * 1.0273863717794265  /2}});
		
		\addplot [RoyalBlue3, line width = 1, smooth, domain=0.6:1.2] {exp{-x*x/2}};
		
		\draw[very thick, dashed, draw=Red3] (0,0.7) -- (1,0.7) node[circle,fill=Red3, inner sep=1.5pt]{} -- (1, {exp{-1.0273863717794265*1.0273863717794265/2}});
		\draw (0.95 ,0.71) node {$\phi_1$};
		\draw (1.015 ,0.67) node {$\psi_1$};
		
		\draw[very thick, dashed, draw=Green3] (0,0.65) -- (0.88,0.65) node[circle,fill=Green3, inner sep=1.5pt]{} -- (0.88, {exp{-1.0273863717794265*1.0273863717794265/2}});
		\draw (0.83 ,0.66) node {$\phi_2$};
		\draw (0.895 ,0.62) node {$\psi_2$};

	\end{axis} 
\end{tikzpicture}


%% file: reflections.tex
\usetikzlibrary{arrows.meta, positioning}

\begin{tikzpicture}[
	circ/.style = {circle, minimum size=1cm, inner sep=0pt},
	arrow/.style = {-Straight Barb, thick, dashed} 
	]
	
	\draw[very thick] (12,1.02) -- (12,-1.02);
	\draw[arrow,draw=LightBlue3]  (0.5,-0.5) -- (9.5,-0.5);
	\draw[arrow,draw=PaleGreen3]  (4.49,0.44) -- (12,0.44);
	\draw[arrow,draw=PaleGreen3] (12,0.56) -- (10.49, 0.56);
	\draw[-Straight Barb, very thick] (-0.5,-1.02) -- (12.5, -1.02);
	
	\node at (12, -1.3) {\large$L$};
	
	\node[circ,fill=LightBlue3] (a) at (0,-0.5) {$X_{p,\tau}$};
	\node[circ,fill=LightBlue3] (b) at (10,-0.5)  {$X_{p,\tau+1}$};
	\node[circ,fill=PaleGreen3] (c) at (4,0.5)   {$X_{p^\prime\!\!,\tau}$};
	\node[circ,fill=PaleGreen3] (d) at (10,0.5)  {$X_{p^\prime\!\!,\tau+1}$};
\end{tikzpicture}